\documentclass[11pt]{article}
\usepackage{rotating}
  \usepackage{graphicx}
  \usepackage{psfrag}
  \usepackage{amsfonts}
\usepackage{amsmath,amssymb}
\usepackage{url}

\usepackage{natbib}
\usepackage{color}

\newtheorem{theorem}{Theorem}[section]

\newcommand{\qed}{\nobreak \ifvmode \relax \else
      \ifdim\lastskip<1.5em \hskip-\lastskip
      \hskip1.5em plus0em minus0.5em \fi \nobreak
      \vrule height0.75em width0.5em depth0.25em\fi}

\newcommand{\R}{\mathrm{I\negthinspace R}}

\newcommand{\st}{\mathbin |}
\newcommand{\bigst}{\mathbin \big |}

\def\comment#1{\textit{[#1]}}
\def\comment#1{}

\begin{document}

\title{Partitioning the Sample Space on Five Taxa for the Neighbor Joining Algorithm}
\author{Kord Eickmeyer\footnote{eMail: kord.eickmeyer@gmx.de} \and Ruriko Yoshida \footnote{eMail: ruriko@ms.uky.edu}}
\date{}
\maketitle
\begin{abstract}
In this paper, we will analyze the behavior of the Neighbor Joining
algorithm on five taxa and we will show that the partition of the
sample (data) space for estimation of a tree topology with five taxa
into subspaces, within each of which the Neighbor Joining algorithm
returns the same tree topology.  A key of our method to partition the
sample space is the action of the symmetric group $S_5$ on the set of
distance matrices by changing the labels of leaves.  The method
described in this paper can be generalized to trees with more than
five taxa.
\end{abstract}

\section{Introduction}
The Neighbor Joining (NJ) algorithm is introduced by \cite{Saitou1987}
and it is widely used because of its accuracy and computational speed.
A number of attempts have been made to understand the good results
obtained with the Neighbor Joining algorithm, especially given the
problems with the inference procedures used for estimating pairwise
distances.  For example, \cite{Bryant2005} showed that the {\em
Q-criterion} (defined in \eqref{eqn:qcrit} in Section
\ref{qcriterion}) is in fact the unique selection criterion such that
it is linear, permutation equivariant, and {\em consistent}, i.e. it
correctly finds the tree corresponding to a tree metric and
\cite{Steel2006} showed a nice review of how the NJ algorithm works.

One of the most important questions to study the behavior of the NJ
algorithm is to analyze its performance with pairwise distances that
are not tree metrics, especially when all pairwise distances are
estimated via the maximum likelihood estimation (MLE).  The following
theorem due to \cite{Atteson99} is one of the results to address this
question:
\begin{theorem}[\cite{Atteson99}]
\label{Atteson}
Neighbor-joining has $l_{\infty}$ radius $\frac{1}{2}$.
\end{theorem}
This means that if the distance estimates are at most half the minimal
edge length of the tree away from their true value then the Neighbor
Joining algorithm will reconstruct the correct tree.  However,
\cite{Levy2005} noted that Atteson's criterion frequently fails to be
satisfied even though the NJ algorithm returns the true tree topology.
Recent work of \cite{Mihaescu2006} extended Atteson's work.
\cite{Mihaescu2006} showed that the NJ algorithm returns the true tree
topology when it works locally for the quartets in the tree.  This
result gives another criterion when then NJ algorithm returns the
correct tree topology and Atteson's theorem is a special case of
Mihaescu et al's theorem.  The results of \cite{Atteson99} and
\cite{Mihaescu2006}, however, depend on an input data and to analyze
the performance of the NJ algorithm for {\em any} input data we need
to study the sample space (the space of all possible distance
matrices).  Therefore, this justifies to study the sample space which
is the subspace of the vector space $\R^{{n \choose 2}}$, where $n$ is
the number of taxa.  In this paper, we will show the partition of the
sample space for estimation of a tree topology with five taxa $\{0, 1,
2, 3, 4\}$ into subspaces, within each of which the NJ algorithm
returns the same tree topology.  For five taxa, the sample space is
the subspace $\R^5$ of the vector space $\R \oplus \R^4 \oplus
\R^5$. Note that our method to partition the sample space into
subspaces described in this paper can be generalized to trees with
more than five taxa.

For every input distance matrix, the NJ algorithm returns a certain
tree topology. It may happen that the minimum Q-criterion is taken by
more than one pair of taxa at some step. In practice, the NJ algorithm
will then have to choose one cherry in order to return a definite
result, but for our analysis we assume that in those cases the NJ
algorithm will return a set containing all tree topologies resulting
from picking optimal cherries.

There are only finitely many tree topologies, and for every topology
$t$ we get a subset $D_t$ of the sample space (input space) such that for all
distance matrices in $D_t$ one possible answer of the NJ algorithm is
$t$. We aim at describing these sets $D_t$ and the relation between
them.

It turns out that we can break up each $D_t$ into easier subsets, by
considering not only the topology returned by the NJ algorithm but
also the order in which cherries are picked. Consider the tree in
figure \ref{fig:fivetree}(a). There are two possible ways in which the
NJ algorithm may reconstruct this tree, namely by first picking the
cherry 1-0 or by first picking the cherry 4-3. The set of input
distance matrices (vectors in $\R \oplus \R^4 \oplus \R^5$) for which
the NJ algorithm will pick exactly the same cherries in the same order
then it turns out to be a cone, so it can easily be described by a set of
halfspaces.

One important tool in our analysis will be the action of the symmetric
group $S_5$ on the set of distance matrices by changing the labels of
leaves. Because there is only one unlabeled tree topology in the case
of five taxa, any two labeled tree topologies can be mapped onto each
other by changing leaf labels. Thus it suffices to analyze just one
$D_t$.

\section{The Neighbor-Joining algorithm}
\label{sec:nj}

\subsection{Input data}
The NJ algorithm is a {\em distance based method} which takes a {\em
distance matrix}, a symmetric matrix $(d_{ij})_{0\leq i, j\leq n-1}$
with $d_{ii} = 0$ representing pairwise distances of a set of $n$ taxa
$\{0, 1, \ldots, n - 1\}$, as the input.  Through this paper, we do
not assume anything on an input data except it is symmetry and $d_{ii}
= 0$. Because of symmetry, the input can be seen as a vector of
dimension $m := \binom{n}{2} = \frac{1}{2}n(n-1)$. We arrange the
entries row-wise, so in the case $n = 4$ we get:
\begin{equation*}
\left(\begin{array}{cccc}
0   & d_0 & d_1 & d_3 \\
d_0 & 0   & d_2 & d_4 \\
d_1 & d_2 & 0   & d_5 \\
d_3 & d_4 & d_5 & 0
\end{array}\right).
\end{equation*}
We denote row/column-indices by pairs of letters such as $a$, $b$, $c$, $d$, while denoting single indices into the ``flattened'' vector by letters $i, j, \dots$. The two indexing methods are used simultaneously in the hope that no confusion will arise. Thus, in the four taxa example we have $d_{0,1} = d_{1,0} = d_0$. In general, we get $d_i = d_{a,b} = d_{b,a}$ with
\begin{equation*}
a = \max\left\{k \bigst \frac{1}{2}k(k-1) \leq i\right\} = 
\left\lfloor \frac{1}{2}+ \sqrt{\frac{1}{4}+2i} \right\rfloor,
b = i-\frac{1}{2}(a-1)a,
\end{equation*}
and for $c>d$ we get
\begin{equation*}
d_{c,d} = d_{c(c-1)/2 + d}.
\end{equation*}

\subsection{The Q-Criterion}\label{qcriterion}

The NJ algorithm starts by computing the so called \emph{Q-criterion} or the {\em cherry picking criterion}, given by the formula
\begin{equation}
\label{eqn:qcrit}\tag{Q}
q_{a,b} := (n-2)d_{a,b} - \sum_{k=0}^{n-1} d_{a,k} - \sum_{k=0}^{n-1} d_{k,b}.
\end{equation}
This is a key of the NJ algorithm to choose which pair of taxa is a neighbor.
\begin{theorem}[\cite{Saitou1987, Studier1988}]
\label{Saitou}
Let $d_{a,b}$ for all pair of taxa $\{a, b\}$ be the tree metric corresponding to the tree $T$.  Then the pair $\{x. y\}$ which minimizes $q_{a,b}$ for all pair of taxa $\{a, b\}$ forms a neighbor.
\end{theorem}

The resulting matrix is again symmetric, and ignoring the diagonal entries we can see it as a vector of dimension $m$ just like the input data. Moreover, the Q-criterion is obtained from the input data by a linear transformation:
\begin{equation*}
\mathbf{q} = A^{(n)}\mathbf{d},
\end{equation*}
and the entries of the matrix $A^{(n)}$ are given by
\begin{equation}
\label{eqn:adef}
A^{(n)}_{ij} = A^{(n)}_{ab,cd} = \begin{cases}
n-4 & \text{if }i = j,\\
-1  & \text{if }i \not = j \text{ and }\{a,b\}\cap \{c,d\} \not = \emptyset,\\
0   & \text{else},
\end{cases}
\end{equation}
where $a > b$ is the row/column-index equivalent to $i$ and likewise for $c > d$ and $j$. When no confusion arises about the number of taxa, we abbreviate $A^{(n)}$ to $A$. In the case of four taxa, we get
\begin{equation*}
A^{(4)} = \left(\begin{array}{cccccc}
0  & -1 & -1 & -1 & -1 & 0\\
-1 & 0  & -1 & -1 & 0  & -1\\
-1 & -1 &  0 &  0 & -1 & -1\\
-1 & -1 &  0 &  0 & -1 & -1\\
-1 & 0  & -1 & -1 & 0  & -1\\
0  & -1 & -1 & -1 & -1 & 0
\end{array}\right).
\end{equation*}

After computing the Q-criterion $\mathbf{q}$, the NJ algorithm
proceeds by finding the minimum entry of $\mathbf{q}$, or, equivalently, the
maximum entry of $-\mathbf{q}$. The two nodes forming the chosen pair
(there may be several pairs with minimal Q-criterion) are then joined
(``cherry picking''), i.e. they are removed from the set of nodes and
a new node is created.  Suppose out of our $n$ taxa
$\{0,\ldots,n-1\}$, the first cherry to be picked is $m-1$, so the
taxa $n-2$ and $n-1$ are joined to form a new node, which we view as
the new node number $n-2$. The reduced pairwise distance matrix is one
row and one column shorter than the original one, and by our choice of
which cherry we picked, only the entries in the rightmost column and
bottom row differ from the original ones. Explicitly,
\begin{equation*}
{d'}_i =
\begin{cases}
d_i & \text{for }0 \leq i < \binom{n-2}{2}
\\
\frac{1}{2}(d_i + d_{i+(n-2)} - d_{m-1})& \text{for }\binom{n-2}{2} \leq i < \binom{n-1}{2} 
\end{cases}
\end{equation*}
and we see that the reduced distance matrix depends linearly on the original one:
\begin{equation*}
\mathbf{d'} = R\mathbf{d},
\end{equation*}
with $R = (r_{ij}) \in \R^{(m-n+1)\times m}$, where
\begin{equation*}
r_{ij} = \begin{cases}
1   & \text{for }0 \leq i=j < \binom{n-2}{2}\\
1/2 & \text{for }\binom{n-2}{2} \leq i < \binom{n-1}{2}, j = i\\
1/2 & \text{for }\binom{n-2}{2} \leq i < \binom{n-1}{2}, j = i+n-2\\
-1/2 & \text{for }\binom{n-2}{2} \leq i < \binom{n-1}{2}, j = m-1\\
0 & \text{else}
\end{cases}
\end{equation*}
The process of picking cherries is repeated until there are only three
taxa left, which are then joined to a single new node.

\section{The cones $C_{ab,c}$}

\begin{figure}[ht]
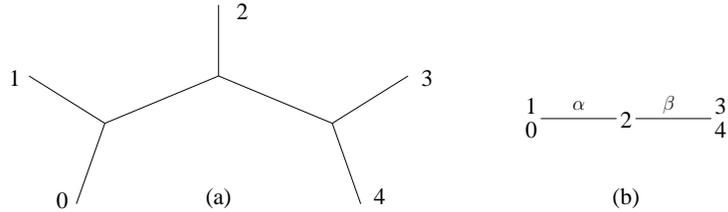

\begin{center}
\resizebox{0.6\textwidth}{!}{\input fivetree.pstex_t}
\end{center}
\caption{(a) A tree with five taxa (b) The same tree with all edges
adjacent to leaves reduced to length zero. The remaining two edges
have lengths $\alpha$ and $\beta$.}
\label{fig:fivetree}
\end{figure}

In the case of five taxa there is just one unlabeled tree topology
(cf. figure \ref{fig:fivetree}) and there are 15 distinct labeled
trees: We have five choices for the leaf which is not part of a cherry
and then three choices how to group the remaining four leaves into two
pairs. For each of these labeled topologies, there are two ways in
which they might be reconstructed by the NJ algorithm: There are two
pairs, any one of which might be chosen in the first step of the NJ
algorithm.

For distinct leaf labels $a$, $b$ and $c \in \{0,1,2,3,4\}$ we define
$C_{ab,c}$ to be the set of all input vectors for which the cherry
$a$-$b$ is picked in the first step and $c$ remains as single node not
part of a cherry after the second step. For example, the tree in
figure \ref{fig:fivetree}(a) is the result for all vectors in
$C_{10,2} \cup C_{43,2}$. Since for each tree topology $((a, b), c, (d, e))$
(this tree topology is written in the Newick format) 
for distinct taxa $a, b, c, d, e \in \{0, 1, 2, 3, 4\}$, the NJ algorithm 
returns the same tree topology with any vector in the
union of two cones $C_{ab, c} \cup C_{de, c}$, there are 30 such cones in 
total, and we call the set of these cones $\mathcal{C}$.

\subsection{The shifting lemma}
\label{sec:shifting}

We first note that there is an $n$-dimensional linear subspace of
$\R^m$ which does not affect the outcome of the NJ algorithm (see
\citep{Mihaescu2006}). For a node $a$ we define its \emph{shift
vector} $\mathbf{s}_a$ by
\begin{equation*}
(\mathbf{s}_a)_{b,c} := \begin{cases}
1 & \text{if }a \in \{b,c\}
\\
0 & \text{else}
\end{cases}
\end{equation*}
which represents a tree where the leaf $a$ has distance 1 from all
other leaves and all other distances are zero. The Q-criterion of any
such vector is $-2$ for all pairs, so adding any linear combination of
shift vectors to an input vector does not change the relative values
of the Q-criteria. Also, regardless of which pair of nodes we join,
the reduced distance matrix of a shift vector is again a shift vector
(of lower dimension), whose Q-criterion will also be constant. Thus,
for any input vector $\mathbf{d}$, the behavior of the NJ algorithm on
$\mathbf{d}$ will be the same as on $\mathbf{d} + \mathbf{s}$ if
$\mathbf{s}$ is any linear combination of shift vectors. We call the
subspace generated by shift vectors $S$.

We note that the difference of any two shift vectors is in the kernel
of $A$, and the sum of all shift vectors is the constant vector with
all entries equal to $n$. If we fix a node $a$ then the set
\begin{equation*}
\{\mathbf{s}_a - \mathbf{s}_b \st b \not = a \}
\end{equation*}
is linearly independent. 

\subsection{Permuting leaf labels}
\label{sec:leafperm}

Because there is only one unlabeled tree topology, we can map any
labeled topology to any other labeled topology by only changing the
labels of the leafs. Such a change of labels also permutes the entries
of the distance matrix: For example, the result of swapping leaves 0
and 1 is the distance vector
\begin{equation*}
\mathbf{d'} := (d_0, d_2, d_1, d_4, d_3, d_5, d_7, d_6, d_8, d_9)^\mathrm{T}.
\end{equation*}
In this way, we get an action of the symmetric group $S_5$ on the input
space $R^{10}$, and the permutation $\sigma \in S_5$ maps the cone
$C_{ab,c}$ linearly to the cone
$C_{\sigma(a)\sigma(b),\sigma(c)}$. Therefore any property of the cone
$C_{ab,c}$ which is preserved by linear transformations must be the
same for all cones in $\mathcal{C}$, and it suffices to determine it
for just one cone.

The action of $S_5$ on $R^{10}$ decomposes into irreducible representations by
\begin{equation*}
\underbrace{\R \oplus \R^4}_{= S} \oplus \underbrace{\R^5}_{=: W},
\end{equation*}
where the first summand is the subspace of all constant vectors and
the second one is the kernel of $A^{(5)}$. The sum of these two
subspaces is exactly the space $S$ generated by the shift vectors. The
third summand, which we call $W$, is the orthogonal complement of $S$
and it is spanned by vectors $w_{ab,cd}$ in $W$ with
\begin{equation*}
(w_{ab,cd})_{xy} := \begin{cases}
1 & \text{if }xy = ab\text{ or }xy = cd
\\
-1 & \text{if }xy = ac\text{ or }xy = bd
\\
0 & \text{else}
\end{cases}
\end{equation*}
where $a$, $b$, $c$ and $d$ are pairwise distinct taxa in $\{0, 1, 2, 3, 4\}$
and $(w_{ab,cd})_{xy}$ is the $x$--$y$th coordinate of the vector $w_{ab,cd}$. 
One linearly independent
subset of this is
\begin{equation*}
w_1 := w_{01,34},\quad
w_2 := w_{12,40},\quad
w_3 := w_{23,01},\quad
w_4 := w_{34,12},\quad
w_5 := w_{40,23}.
\end{equation*}
Note that the 5-cycle $(01234)$ of leaf labels cyclically permutes
these basis vectors, whereas the transposition $(01)$ acts via the
matrix
\begin{equation*}
T := \frac{1}{2}
\left(\begin{array}{rrrrr}
2 & 1 & 1 & 1 & 1 \\
0 & 1 &-1 &-1 &-1 \\
0 &-1 &-1 & 1 &-1 \\
0 &-1 & 1 &-1 &-1 \\
0 &-1 &-1 &-1 & 1
\end{array}\right).
\end{equation*}
Because a five-cycle and a transposition generate $S_5$, in principle
this gives us complete information about the operation.

\subsection{The cone $C_{43,2}$}

Since we can apply a permutation $\sigma \in S_5$, without loss of generality,
we suppose that the first cherry to be picked is the cherry 9, which is the 
cherry with leaves 3 and 4. This is true for all input vectors 
$\mathbf{d}$ which satisfy
\begin{equation*}
(\mathbf{h}_{9,i}, \mathbf{d}) \geq 0
\text{ for }i=0,\ldots,8,
\end{equation*}
where the vector
\begin{equation*}
\mathbf{h}^{(n)}_{ij} := -A^{(n)}(\mathbf{e}_i - \mathbf{e}_j)
\end{equation*}
is perpendicular to the hyperplane of input vector for which cherries
$i$ and $j$ have the same Q-criterion, pointing into the direction of
vectors for which the Q-criterion of cherry $i$ is lower.

We let $\mathbf{r}_1$, $\mathbf{r}_2$ and $\mathbf{r}_3$ be the first
three rows of $-A^{(4)}R^{(5)}$. If $(\mathbf{r}_1, \mathbf{d})$ is
maximal then the second cherry to be picked is 0-1, leaving 2 as the
non-cherry node, and similarly $\mathbf{r}_2$ and $\mathbf{r}_3$ lead
to non-cherry nodes 1 and 0. This allows us to define the set of all
input vectors $\mathbf{d}$ for which the first picked cherry is 3-4
and the second one is 0-1:
\begin{equation}
\label{eqn:top34-2}
C_{34,2} := \{ \mathbf{d} \st
(\mathbf{h}_{9,i}, \mathbf{d}) \geq 0
\text{ for }i=0,\ldots,8,
\text{ and }
(\mathbf{r}_3-\mathbf{r}_1,\mathbf{d}) \geq 0,
(\mathbf{r}_3-\mathbf{r}_2,\mathbf{d}) \geq 0
\}.
\end{equation}

We have defined this set by 11 bounding hyperplanes. However, in fact, the
resulting cone has only nine facets. A computation using
\texttt{polymake} \citep{Gawrilow2000} reveals that the two hyperplanes 
$\mathbf{h}_{9,1}$
and $\mathbf{h}_{9,2}$ are no longer faces of the cone, while the
other nine hyperplanes in \eqref{eqn:top34-2} give exactly the facets
of the cone.  

\subsection{The Rays of $\mathcal{C}$}

Again using \texttt{polymake} \citep{Gawrilow2000}
we find that $C_{43,2}$ is the positive
cone spanned by fourteen rays. The union of the orbits of these rays
under the action of $S_5$ is a set of 82 vertices, which we call
$\mathcal{R}$. Each of the cones in $\mathcal{C}$ is spanned by a
certain subset with 14 vertices each of $\mathcal{R}$.

\begin{figure}[ht]
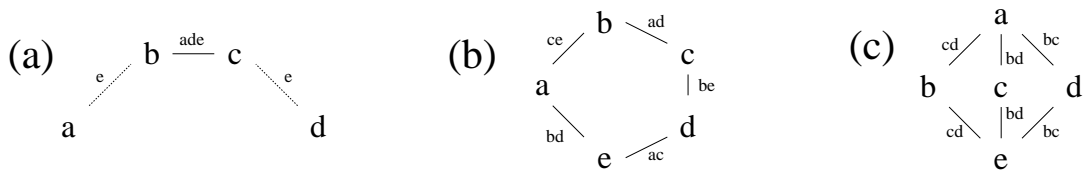

\begin{center}
\resizebox{0.9\textwidth}{!}{\input rays.pstex_t}
\end{center}
\caption{Graphs characterizing the rays in $\mathcal{R}$. Here $a$,
$b$, $c$, $d$ and $e$ are pairwise distinct labels between 0 and 4.}
\label{fig:rays}
\end{figure}

The set $\mathcal{R}$ has three orbits under the action of the
symmetric group $S_5$.  We characterize the rays in these orbits using the
graphs in Figure \ref{fig:rays}. In these graphs, nodes which are
connected by an edge form a pair with the minimum Q-criterion in the first
step of the NJ algorithm, and the labels to each edge show which nodes 
are possible as the remaining unpaired node.

For each graph in figure \ref{fig:rays} and each assignment of the
five leaves $0,\ldots,4$ to the variables $a,\ldots,e$ we get a set
$G$ of vectors which belong exactly to the cones indicated by the
graph. By what we have said in sections \ref{sec:shifting} and
\ref{sec:leafperm}, the set $G$ has the form
\begin{equation*}
  G = \{ \alpha\mathbf{g} + s \st \alpha \geq 0\text{ and }s\in S\},
\end{equation*}
where $S$ is the set of all shifting vectors in the sample space (which
is the subspace in $R^{10}$).
Thus $G$ is described by a vector $\mathbf{g}$ which is unique up to
normalization by a positive constant. We call this vector
$\mathbf{g}^{(\mathrm{a})}_{abcde}$,
$\mathbf{g}^{(\mathrm{b})}_{abcde}$ or
$\mathbf{g}^{(\mathrm{c})}_{abcde}$, depending on which graph it
corresponds to.

Some of these graphs coincide for different choices of
$a,\ldots,e$.
\begin{itemize}
\item The graphs in \ref{fig:rays}(a) are strings of four taxa. There
are five ways to choose the four taxa, and $\frac{4!}{2}=12$ ways in
which to arrange them (reversing the order of the string does not
matter). In total there are 60 rays of this type.

Thus for example the string 0-1-2-3 refers to an input vector
$\mathbf{g}^{(\mathrm{a})}_{01234}$ which
lies in the intersection of the five cones $C_{10,4}, C_{21,0}, C_{21,3}, 
C_{21,4}, C_{30, 4}$ but not in any other cone. 
It is worth noting that this vector is a
ray of the cones $C_{21,0}, C_{21,3}$ and $C_{21,4}$ but \emph{not} of 
the cones $C_{10,4}$
and $C_{32,4}$, which is indicated by the dotted edges in figure
\ref{fig:rays}. This implies that the cones in $\mathcal{C}$ do
\emph{not} form a fan, and in particular no polytope in $\R^{10}$ can have
$\mathcal{C}$ as its normal fan.

\item There are $\frac{4!}{2}$ possibilities for the cycle graph in
figure \ref{fig:rays}(b), each of which represents an input vector
which forms a ray of ten of the cones in $\mathcal{C}$.

\item Finally, there are $\binom{5}{2}=10$ graphs of type
\ref{fig:rays}(c), because these are determined by the two nodes we
choose for $a$ and $e$. Each of these corresponds to an input vector
that forms a ray of twelve cones.
\end{itemize}
If we collect all those rays in $\mathcal{R}$ which are rays of
$C_{43,2}$ we get
\begin{itemize}
\item 6 rays of type (a):
\begin{equation*}
\mathbf{g}^{(\mathrm{a})}_{03412},
\mathbf{g}^{(\mathrm{a})}_{03421},
\mathbf{g}^{(\mathrm{a})}_{13402},
\mathbf{g}^{(\mathrm{a})}_{13420},
\mathbf{g}^{(\mathrm{a})}_{23401},
\mathbf{g}^{(\mathrm{a})}_{23410}
\end{equation*}

\item 4 of type (b):
\begin{equation*}
\mathbf{g}^{(\mathrm{b})}_{34201},
\mathbf{g}^{(\mathrm{b})}_{34210},
\mathbf{g}^{(\mathrm{b})}_{34012},
\mathbf{g}^{(\mathrm{b})}_{34102}
\end{equation*}

\item 4 of type (c):
\begin{equation*}
\mathbf{g}^{(\mathrm{c})}_{34201},
\mathbf{g}^{(\mathrm{c})}_{34210},
\mathbf{g}^{(\mathrm{c})}_{34021},
\mathbf{g}^{(\mathrm{c})}_{34120}.
\end{equation*}
\end{itemize}
We call these vectors $\mathbf{g}_1,\ldots,\mathbf{g}_{14}$.
Any vector $\mathbf{v} \in C_{43,2}$ is a linear combination of these
fourteen vectors with nonnegative linear coefficients, say
\begin{equation*}
\mathbf{v} = \sum_i \alpha_i \mathbf{g}_i,
\end{equation*}
where we number the fourteen rays arbitrarily. By linearity, in both
steps of the NJ algorithm, the Q-criteria are linear combinations of 
the Q-criteria
for the $\mathbf{g}_i$, using the same linear coefficients
$\alpha_i$. Furthermore, because all of the vectors $\mathbf{g}_i$ lie
in $C_{43,2}$, in both steps of the NJ algorithm there is a pair of leaves 
such that
the minimum of the Q-criteria is attained at this pair for all
$\mathbf{g}_i$, and therefore the minimum Q-criterion of the linear
combination $\sum \alpha_i \mathbf{g}_i$ is the linear combination of
the minimum Q-criteria of the individual $\mathbf{g}_i$. This minimum
is attained exactly by those pairs which have minimum Q-criterion in all
$\mathbf{g}_i$ with strictly positive linear coefficient $\alpha_i$,
we get:
\begin{equation*}
\mathbf{v} \in C_{ab,c}
\quad
\text{iff}
\quad
\mathbf{g}_i \in C_{ab,c}\text{ for all }i\text{ with }\alpha_i > 0.
\end{equation*}

\begin{figure}[ht]
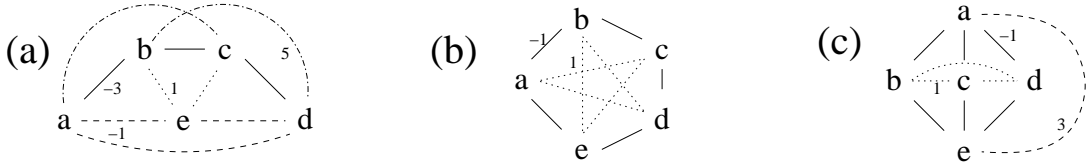

\begin{center}
\resizebox{0.9\textwidth}{!}{\input rayvectors.pstex_t}
\end{center}
\caption{Graphs describing the vectors $\mathbf{g}_{abcde}$ for each
of the three types of graphs. Labels assigned to edges denote the
distance assigned to the corresponding pair of leaves in the vector,
and edges of the same type receive the same label.}
\label{fig:rayvectors}
\end{figure}

Therefore the rays in $\mathcal{R}$ are the input vectors for which the NJ
algorithm is least stable. We can give explicit descriptions of these vectors
using the graphs in Figure \ref{fig:rayvectors}. The graphs give the
distance value assigned to each pair of leaves. For example, we get
\begin{equation*}
\mathbf{g}^{(\mathrm{a})}_{01234} = (-3, 5, -3, -1, 5, -3, -1, 1, 1, -1 )^\mathrm{T}.
\end{equation*}

\subsection{The Cone of Tree Metrics}

Suppose our input vector is a tree
metric. By relabeling the tree and using the shifting lemma, we can
reduce the tree to Figure \ref{fig:fivetree}(b), where all edges
connecting to leaves have length zero. The edge-lengths $\alpha$ and
$\beta$ determine the distance vector $\mathbf{v}$ as follows:
\begin{equation*}
\mathbf{v} = \alpha(0, 1, 1, 1, 1, 0, 1, 1, 0, 0)^\mathrm{T}
+
\beta(0, 0, 0, 1, 1, 1, 1, 1, 1, 0)^\mathrm{T}
=: \alpha\mathbf{a} + \beta\mathbf{b}.
\end{equation*}
If we assume these to be non-negative, we get a cone spanned by the
two vectors $\mathbf{a}$ and $\mathbf{b}$. These satisfy
\begin{equation*}
\mathbf{a} \in C_{10,2} \cap C_{10,3} \cap C_{10,4},
\quad
\mathbf{b} \in C_{43,0} \cap C_{43,1} \cap C_{43,2}
\quad\text{and}\quad
\underbrace{\frac{1}{2}(\mathbf{a}+\mathbf{b})}_{=: \mathbf{b}'} \in C_{10,2} \cap C_{43,2},
\end{equation*}
and none of these vectors lies in any other cone.

We normalize to $\alpha+\beta = 1$, and because of symmetry we may
assume that $\alpha \geq \beta$, or equivalently $\alpha \geq
\frac{1}{2}$. Then we can write $\mathbf{v}$ as
\begin{align*}
\mathbf{v} 
&= \alpha \mathbf{a} + (1-\alpha)\mathbf{b}
\\
&= \left(\alpha-\frac{1}{2}\right) \mathbf{a} + (1-\alpha)\mathbf{b}'
\end{align*}
Because both $\mathbf{a}$ and $\mathbf{b}'$ are in the convex set
$C_{10,2}$, the convex linear combination $\mathbf{v}$ is also in
$C_{10,2}$, proving correctness of the NJ algorithm for five taxa.

We can also compute the $\ell_2$ distance of $\mathbf{v}$ from the
faces of the cone $C_{10,2}$. There are 9 hyperplanes defining $C_{10,2}$, 
but we may ignore one of them, which defines the common facet with
$C_{43,2}$. For any hyperplane $H$ such that $\mathbf{a}$ and
$\mathbf{b}'$ lie on the same side of $H$, the distance
$d(\mathbf{v},H)$ between $H$ and $\mathbf{v}$ is given by
\begin{equation*}
d(\mathbf{v},H) = \left(\alpha-\frac{1}{2}\right) d(\mathbf{a},H) + (1-\alpha)d(\mathbf{b}',H),
\end{equation*}
and taking the minimum of this over the eight remaining faces of
$C_{10,2}$ we obtain
\begin{equation*}
d(\mathbf{v},(C_{10,2}\cap C_{43,2})^\mathrm{c}) = \frac{1-\alpha}{\sqrt{3}},
\end{equation*}
where $(C_{10,2}\cap C_{43,2})^\mathrm{c}$ is the complement of 
$C_{10,2}\cap C_{43,2}$.
If we divide this by the length $\beta = 1-\alpha$ of the smaller of
the two interior edges, we get an $\ell_2$-radius of $1/\sqrt{3}
\approx 0.577$. Note that because our method relies on orthogonal
projections, we get $\ell_2$ bounds instead of $\ell_\infty$ bounds.


\section{Simulation results}
\begin{figure}[ht]
\begin{center}
\includegraphics[scale=0.6]{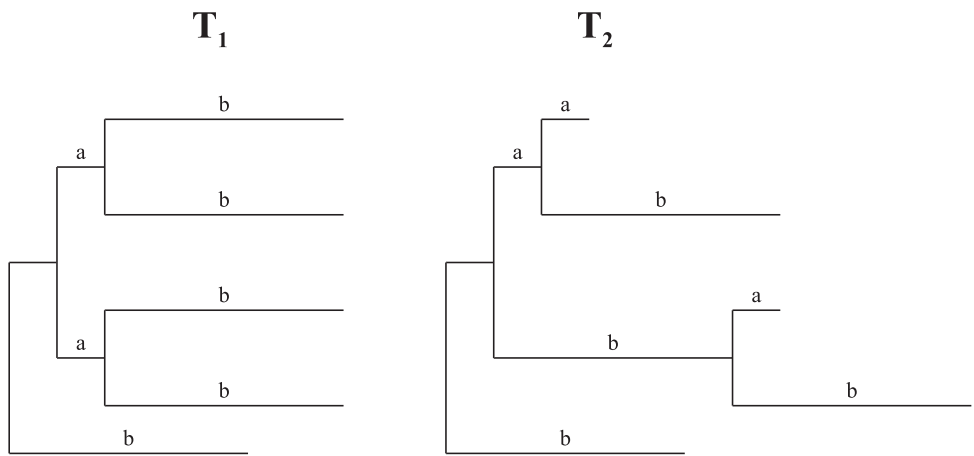}
\end{center}
\caption{$T_1$ and $T_2$ tree models which are subtrees of the tree
models in \citep{Ota2000}.}
\label{fig:T1T2}
\end{figure}

In this section we will analyze how the tree metric for a tree and
pairwise distances estimated via the maximum likelihood estimation
locate in the partition of the sample space.  Particularly, we analyze
subtrees of the two parameter family of trees described by
\citep{Ota2000}.  These are trees for which the NJ algorithm has
difficulty in resolving the correct topology.  In order to understand
how they locate to each other, we simulated 10,000 data sets on each
of the two tree shapes, $T_1$ and $T_2$ (Figure \ref{fig:T1T2}) at the
edge length ratio, a/b = 0.03/0.42 for sequences of length 500BP under
the Jukes-Cantor model \citep{Jukes1969}. We also repeated the runs
with the Kimura 2-parameter model \citep{Kimura1980}. They are the
cases (on eight taxa) in \citep{Ota2000} that the NJ algorithm had
most difficulties in their simulation study (also the same as in
\citep{Levy2005}).  Each set of 5 sequences are generated via {\tt
evolver} from {\tt PAML} package \citep{Yang1997} under the given
model.  {\tt evolver} generates a set of sequences given the model and
tree topology using the birth-and-death process.  For each set of 5
sequences, we compute first pairwise distances via the heuristic MLE
method using a software {\tt fastDNAml} \citep{Olsen1994}.
To compute cones, we used {\tt MAPLE} and \texttt{polymake}.

To study how far each set of pairwise distances estimated via the
maximum likelihood estimation (which is a vector $\mathbf{y}$ in
$\R^5$) locates from the cone, where the additive tree metric lies, in
the sample space, we calculated the $\ell_2$-distance between the cone
and a vector $\mathbf{y}$.

Suppose we have a cone $C$ defined by hyperplanes
$\mathbf{n}_1,\ldots,\mathbf{n}_r$, i.e.
\begin{equation*}
C = \{ \mathbf{x} \st (\mathbf{n}_i,\mathbf{x}) \geq 0\text{ for
}i=1,\ldots,r\},
\end{equation*}
and we want to find a closest point in $C$ from some given point
$\mathbf{v}$. Because $C$ is convex, for $\ell_2$-norm there is only
one such point, which we call $\mathbf{u}$. If $\mathbf{v} \in C$ then
$\mathbf{u} = \mathbf{v}$ and we are done. If not, there is at least
one $\mathbf{n}_i$ with $(\mathbf{n}_i, \mathbf{v}) < 0$, and
$\mathbf{u}$ must satisfy $(\mathbf{n}_i,\mathbf{u}) = 0$.

Now the problem reduces to a lower dimensional problem of the same
kind: We project $\mathbf{v}$ orthogonally into the hyperplane $H$
defined by $(\mathbf{n}_i,\mathbf{x}) = 0$ and call the new vector
$\tilde{\mathbf{v}}$. Also, $C \cap H$ is a facet of $C$, and in
particular a cone, so proceed by finding the closest point in this
cone from $\tilde{\mathbf{v}}$.

\begin{figure}[ht]
\resizebox{0.48\textwidth}{!}{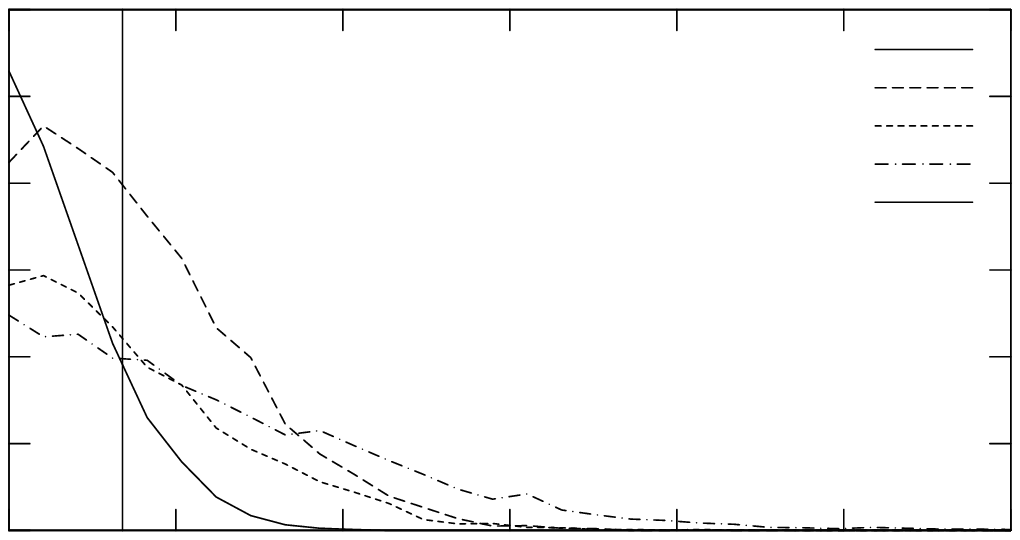}
\resizebox{0.48\textwidth}{!}{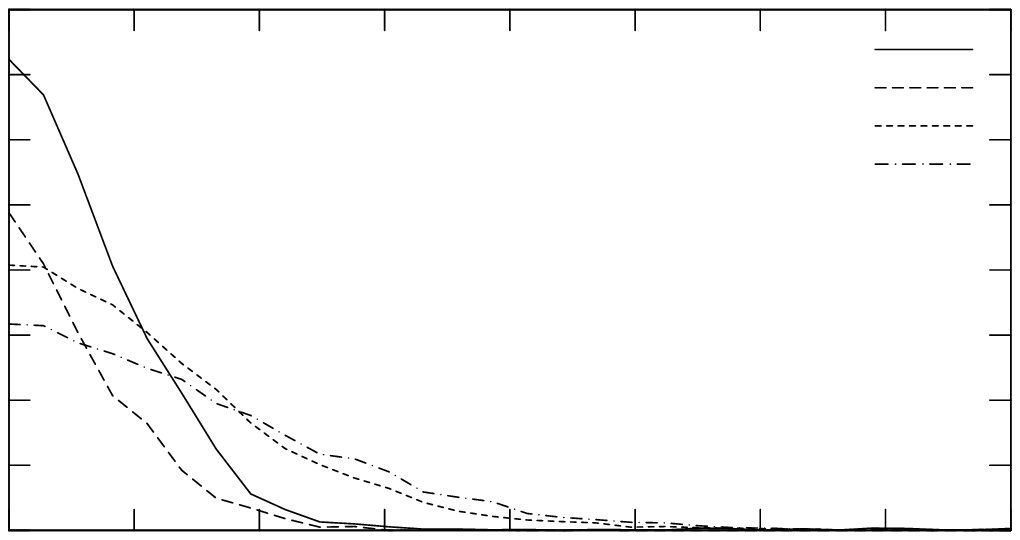}
\caption{Distances of correctly (left) and incorrectly (right)
classified input vectors from the closest incorrectly/correctly
classified vector.}
\label{fig:simplots}
\end{figure}

We say an input vector (distance matrix) is {\em correctly classified}
if the vector locates in one of the cones where the vector representation 
of the tree metric (noiseless input) lies.
We say an input vector is {\em incorrectly classified}
if the vector locates in the complement of the cones where 
the vector representation of the tree metric lies.
For input vectors (distance matrices) which are correctly classified by 
the NJ algorithm, we compute the
minimum distance to any cone giving a different tree topology. This
distance gives a measure of robustness or confidence in the result,
with bigger distances meaning greater reliability. The results are
plotted in the left half of Figure \ref{fig:simplots} and in Figure
\ref{fig:goodstat}. Note that the distance of the noiseless input,
i.e. the tree metric from the tree we used for
generating the data samples, gives an indication of what order of
magnitude to expect with these values.

\begin{figure}[ht]
\begin{center}
\begin{tabular}{l|cc|cc}
&\multicolumn{2}{c}{\textbf{JC}}&\multicolumn{2}{c}{\textbf{Kimura2}}
\\
&T1&T2&T1&T2
\\
\hline
\textbf{\# of cases}& 3,581 & 6,441 & 3,795 & 4,467
\\
\textbf{Mean}&$0.0221$&$0.0421$ & $0.0415$ & $0.0629$
\\
\textbf{Variance}& $2.996\cdot 10^{-4}$ & $9.032\cdot 10^{-4} $
&
$1.034 \cdot 10^{-3}$ & $2.471 \cdot 10^{-3}$
\end{tabular}
\end{center}
\caption{Mean and variance of the
distances of correctly classified
vectors from the nearest
misclassified vector.}
\label{fig:goodstat}
\end{figure}

For input vectors to which the NJ algorithm answers with a tree topology 
different
from the correct tree topology, we compute the
distances to the two cones for which the correct answer is given and
take the minimum of the two. The bigger this distance is, the further
we are off. The results are shown in the right half of Figure
\ref{fig:simplots} and in Figure \ref{fig:badstat}.

\begin{figure}[ht]
\begin{center}
\begin{tabular}{l|cc|cc}
&\multicolumn{2}{c}{\textbf{JC}}&\multicolumn{2}{c}{\textbf{Kimura2}}
\\
&T1&T2&T1&T2
\\
\hline
\textbf{\# of cases}& 6,419 & 3,559 & 6,205 & 5,533
\\
\textbf{Mean}&$0.0594$&$0.0331$ & $0.0951$ & $0.0761$
\\
\textbf{Variance}& $0.0203$ & $7.39\cdot 10^{-4} $
&
$0.0411$ & $3.481 \cdot 10^{-3}$
\end{tabular}
\end{center}
\caption{Mean and variance of the
distances of misclassified vectors to the nearest correctly classified vector.}
\label{fig:badstat}
\end{figure}

From our results in Figure \ref{fig:goodstat} and Figure \ref{fig:badstat},
one notices that the NJ algorithm returns the correct tree more often with 
$T_2$ than with $T_1$.  These results are consistent with the results in 
\citep{Steel2006, Mihaescu2006}.  Note that any possible quartet in $T_1$
has a smaller (or equal) length of its internal edge than in $T_2$ 
(see Figure \ref{fig:T1T2}).  \cite{Steel2006} defined this measure 
as {\em neighborliness}.  \cite{Mihaescu2006} showed that the NJ algorithm
returns the correct tree if it works correctly locally for the quartets in
the tree.  The neighborliness of a quartet is one of the most important
factors to reconstruct the quartet correctly, i.e., the shorter it is the more 
difficult the NJ algorithm returns the correct quartet. 
Also Figure \ref{fig:simplots} shows that most of the input vectors lie around
the border lines of cones, including the noiseless input vector (the tree 
metric).  This shows that the tree models $T_1$ and $T_2$ are difficult
for the NJ algorithms to reconstruct the correct trees.  
All source codes for these simulations described in this paper 
will be available at authors' websites.

\pagebreak
\setcounter{page}{1}
\bibliographystyle{named}
\bibliography{phylo}

\end{document}